\documentclass{elsart}
\usepackage{graphics}
\usepackage{graphicx}
\usepackage{epsfig}
\usepackage{amssymb,amsmath}
\usepackage[colorlinks]{hyperref}

\newtheorem{Lemma}{Lemma}[section]
\newtheorem{Corollary}{Corollary}[section]
\newtheorem{Theorem}{Theorem}[section]
\newtheorem{Proposition}{Proposition}[section]
\newtheorem{Definition}{Definition}[section]
\newtheorem{Remark}{Remark}[section]

\begin{document}
\linespread{1}
\begin{frontmatter}
\title{The Farey map exploited in wavelet analysis - Construction of an associated Farey mother wavelet}
\author{Sabrine Arfaoui\thanksref{label2}}
\address{Laboratory of Algebra, Number Theory and Nonlinear Analysis, LR18ES15, Department of Mathematics, Faculty of Sciences, 5000 Monastir, Tunisia.}
\ead{sabrine.arfaoui@issatm.rnu.tn} 
\author{Riadh Chteoui\thanksref{label2}}
\address{Laboratory of Algebra, Number Theory and Nonlinear Analysis, LR18ES15, Department of Mathematics, Faculty of Sciences, 5000 Monastir, Tunisia.}
\ead{riadh.chteoui.fsm@gmail.com}
\author{Anouar Ben Mabrouk\corauthref{cor1}\thanksref{label1}\thanksref{label2}}
\address{Department of Mathematics, Higher Institute of Applied Mathematics and Computer Science, University of Kairoaun, Street of Assad Ibn Alfourat, 3100 Kairouan, Tunisia.}
\ead{anouar.benmabrouk@fsm.rnu.tn}
\thanks[label1]{Laboratory of Algebra, Number Theory and Nonlinear Analysis, LR18ES15, Department of Mathematics, Faculty of Sciences, 5000 Monastir, Tunisia.}
\thanks[label2]{{Department of Mathematics, Faculty of Sciences, University of Tabuk, Saudi Arabia.}}
\corauth[cor1]{Corresponding author.}
\begin{abstract}
In the present work the well known Farey map is exploited to consruct a new mother wavelet. Properties such as admissibility, moments, 2-scale relation and reconstruction rule have been established. The constructed mother may be a good candidate to analyze hyperbolic problems such as hyperbolic PDEs.
\end{abstract}
\end{frontmatter}
\section{Introduction}
Wavelet analysis was introduced since the 80th decade of the last century in petroleum exploration as a refinement of Fourier analysis which failed in extracting the characteristics of signals subjects of the exploration. Since then, wavelets have attracted the interest of researchers in both pure and applied fields.

For a large community, of non mathematicians, a wavelet may be defined as a wave function which decays rapidly and which has besides a zero mean. Wavelet analysis consists of breaking up a signal into parts relatively to approximating functions obtained as shifted and dilated versions of the wavelet \cite{Daubechies}.

In Fourier analysis, description of signals is restricted to the global behavior and can not provide information about hidden details due to its lack of time-frequency and/or time-space localization. In 1940, Denis Gabor introduced the so-called windowed Fourier transform (WFT) to overcome this lack. The Fourier modes used in Fourier transform are multiplied by suitable functions localized in time such as Gaussian window. This step permitted to understand new situations. But the situation has changed again with the discovery and/or the emergence of new problems related to irregular variations such as gravitational waves where glitches (which are are bursts of noise) remain after filtering. This leads researchers to think about more sophisticated tools for signal processing and leads next to the discovery and development of wavelet analysis. Wavelets thus permit the localization of analyzed signals in both time and frequency. Contrarily to Fourier analysis, wavelets permit also to analyze non-stationarity, non seasonality, irregularity with more precision.

In the present work we propose to provide a rigorous development of new wavelet mother by exploiting the characteristics of the well known Farey map. The remaining parts of the present document will be organized as follows. Section 2 is devoted to the review of wavelet analysis. Section is concerned with the development of our main results dealing with a new type of wavelet function based on the well-known Farey map which we call Farey wavelet. Special characteristics such as admissibility, moments, reconstruction rule have been established. Section 4 is a conclusion.
\section{Wavelet revisited}
A wavelet is a wave function that decays rapidly and has a zero average value. Wavelet analysis is a breaking up of a signal into approximating functions (shifted and dilated versions of the wavelet) contained in finite domains \cite{Daubechies}.
The wavelet analysis was introduced in the early 1980s in the context of signal analysis and petroleum exploration. It aims to give a representation of signals and detect their characteristic. Several methods previously have been used, the most known is the Fourier transform. In Fourier transform further description of signals is limited to the overall behavior and can not provide any information on the details. In digital signal processing, Fourier analysis often requires linear calculation algorithms. In 1940, Denis Gabor introduced the windowed Fourier transform (TFF) to address the problems of time-frequency localization. It consists in calculating the Fourier transform of the signal by multiplying it by a function localized in time (Gaussian window) and then calculating the transform. But the situation has changed with the emergence of new problems especially irregular variation. The major drawback of the TFF is the shape stability and the window's size. Gaussian type windows can not for example model non stationary properties. Henceforth the need for analysis using non-linear algorithms, non-stationary signals and/or non-periodical bases became necessary. Specifically analysis with well-specified characteristics, i.e. localization time and frequency, adaptivity to the data, easily implemented advanced algorithms and optimum computation time needs to be developed. This was how wavelet theory was born. It had subsequently renewed interest and has been steadily developed in theory and application. Wavelets differ from Fourier methods in that they allow the localization of a signal in both time and frequency. It is a tool which breaks up data into different frequency components or sub bands and then studies each component with a resolution that is matched to specific or proper scale. Unlike the Fourier series, it can be used on non-stationary transient signals with more precise results \cite{Othman}. This section is devoted to present the main ideas on wavelet analysis namely wavelet transforms, multi-resolution analysis, wavelet bases and algorithms of construction and reconstruction.

In purely mathematical point of view, a wavelet is a function $\psi\in\,L^2(\mathbb{R})$  which satisfies the following conditions.
\begin{itemize}
	\item Admissibility,
	\begin{equation}\label{eqn:ond}
\mathcal{A}_\psi=\displaystyle\int_{\mathbb{R}^+}|\hat\psi(\omega)|^2{d\omega\over{|\omega|}}<\infty.
	\end{equation}
	\item Zero mean,
	\begin{equation}
	\widehat{\psi}(0)=\displaystyle\int_{-\infty}^{+\infty}\psi(u)du=0.
	\end{equation}
	\item Localization,
	\begin{equation}
	\|\psi\|_2^2=\displaystyle\int_{-\infty}^{+\infty}|\psi(u)|^{2}du=1.
	\end{equation}
	\item Enough vanishing moments,
	\begin{equation}
	p=0,...,m-1,\quad\displaystyle\int_{\mathbb{R}}\psi(t)t^pdt=0.
	\end{equation}
\end{itemize}
To analyze functions/signals by wavelets, one passes by the so-called wavelet transforms. A wavelet transform (WT) is a representation of a time-frequency signal. It replaces Fourier sine/cosine by a wavelet. Generally, there are two types of processing; The continuous wavelet transform and discrete wavelet transform.

The CWT is based firstly on the introduction of a translation parameter $u\in\mathbb{R}$ and another parameter $s>0$ known as the scale to the analyzing wavelet $\psi$ called the mother wavelet. The translation parameter determines the position or the time around which we want to assess the behavior of the signal, while the scale factor is used to assess the signal behavior around the position. That is, it allows us to estimate the frequency of the signal at that point. Let
\begin{equation}
\psi_{s,u}(x)=\displaystyle\frac{1}{\sqrt{s}}\psi(\displaystyle\frac{x-u}{s}).
\end{equation}
The continuous wavelet transform at the position $u$ and the scale $s$ is defined by
\begin{equation}\label{CWT}
d_{u,s}(f)=\displaystyle\int_{-\infty}^{\infty}\psi_{u,s}(t)f(t)dt,\quad\forall\,u,s.
\end{equation}
By varying the parameters $s$ and $u$, we may cover completely all the time-frequency plan. This gives a full and redundant representation of the whole signal to be analyzed (See Mallat (2000) \cite{Mallat2000}). This transform is called continuous because of the nature of the parameters $s$ and $u$ that may operate at all levels and positions. The original function $f$ can be reproduced knowing its CWT by the following relationship.
\begin{equation}\label{reconstructionformula}
f(x)=\dfrac{1}{\mathcal{A}_\psi}\displaystyle\int\displaystyle\int_{\mathbb{R}}d_{u,s}(f)\psi\left(\displaystyle\frac{x-u}{s}\right)\displaystyle\frac{dsdu}{s^2}.
\end{equation}
The DWT is obtained by restricting the scale qnd position parameters to a discrete grid. The most known method is the dyadic grid $s=2^{-j}$ and $u=k2^{-j}$, $j,k\in\mathbb{Z}$. In this case, the wavelet copy $\psi_{u,s}$ is usually denoted by $\psi_{j,k}(x)=2^{-j/2}\psi(2^jx-k)$. The DWT of a function $f$ is 
\begin{equation}\label{DWT}
d_{j,k}=\displaystyle\int_{-\infty}^{\infty}\psi_{j,k}(t)S(t)dt.
\end{equation}
These are often called wavelet coefficients or detail coefficients of the signal $S$.

It holds that the set $\left(\psi_{j,k}\right)_{j,k\in\mathbb{Z}}$ constitutes an orthonormal basis of $L^2(\mathbb{R})$ and called wavelet basis. An element $f\in L^2(\mathbb{R})$ is decomposed according to this basis into a series
\begin{equation}\label{waveletseries}
f(t)=\displaystyle\sum_{j=0}^{\infty}\displaystyle\sum_{k}d_{j,k}\psi_{j,k}(t)
\end{equation}
called the wavelet series of $f$ which replaces the reconstruction formula for the CWT.

It holds in wavelet theory that the previous concepts may induce a functional framework for representing functions by a series of approximations called resolutions. Such framework is known as the multiresolution analysis (MRA) on $\mathbb{R}$. MRA is a family of closed vector subspaces $(V_j)_{j\in\mathbb{Z}}$ of $L^2(\mathbb{R})$. For each $j\in\mathbb{Z}$, $V_j$ is called the approximation at the scale or the level $j$. More precisely (\cite{Mallat2000}),	a multi-resolution analysis is a countable set of closed subsets $(V_j)_{j\in\mathbb{Z}}$ of $L^2(\mathbb{R})$ that satisfies the following points.
	\begin{itemize}
		\item[a)]$\forall j\in\mathbb{Z}$; $...V_{-j-1}\subset V_{-j}\subset...\subset V_{-1}\subset V_{0}\subset V_{1}\subset....\subset V_{ j}\subset V_{j+1}$.
		\item[b)] $\displaystyle\bigcap_{j\in\mathbb{Z}} V_j = \{0\}$.
		\item[c)] $\displaystyle\overline{\bigcup_{j\in\mathbb{Z}}V_j}=L^2(\mathbb{R})$.
		\item[d)] $\forall j\in \mathbb{Z}$; $f\in V_j\Leftrightarrow f(2.)\in V_{j+1}$
		\item[e)] $\forall j\in \mathbb{Z}$; $f\in V_j\Leftrightarrow f(x-k)\in V_{j}$
		\item[f)] There exists $\varphi\in\,V_0$ such that $\left\{\varphi_{0,k}=\varphi(.-k);\,k\in\mathbb{Z}\right\}$ is an orthogonal Riesz basis of $V_0$.
	\end{itemize}
The source function $\varphi$ is called the scaling function of the MRA or also the wavelet father.

It holds that this function generates all the subspaces $V_j$'s of the MRA by acting dilation/translation parameters. Indeed, the set $\left(\varphi_{j,k}(x)=2^{j/2}\varphi(2^{j}x-k)\right)_{k}$ is an orthogonal basis of $V_j$ for all $j,\in\mathbb{Z}$. Moreover, there is an orthogonal supplementary $W_j$ of $V_j$ in $V_{j+1}$, that is
\begin{equation}\label{Vj+1VjWj}
V_{j+1}=V_j\oplus^\bot W_j
\end{equation}
The space $W_j$, $j\in\mathbb{Z}$ is called detail space at the scale or the level $j$ for which the set $\left(\psi_{j,k}(t)=2^{j/2}\psi(2^{j}t-k)\right)_{k}$ is an orthogonal basis.

The strongest point in MRA and wavelet theory is that the scaling function and the analyzing wavelet leads each one to the other. Indeed, recall that $\varphi$ belongs to $V_0\subset\,V_1$ and the latter is generated by the basis $(\varphi_{1,k})_k$. Hence, $\varphi$ is expressed by means of $(\varphi_{1,k})_k$. More precisely, we have the so-called 2-scale relation
\begin{equation}\label{equation2echelle}
\varphi(x)=\sqrt{2}\displaystyle\sum_{k}h_{k}\varphi(2x-k)
\end{equation}
where the coefficients $h_k$ are
$$
h_{k}=\displaystyle\int_{\mathbb{R}}\varphi(x)\overline{\varphi(2x-k)}dx.
$$
It holds that the mother wavelet $\psi$ satisfies 
$$
\psi(x)=\sqrt{2}\displaystyle\sum_{k}g_{k}\varphi(2x-k)
$$
where the $g_k$'s are evaluated by
$$
g_k=(-1)^kh_{1-k}.
$$
For more details, we refer to \cite{Daubechies}, \cite{Hardleetal1997}, \cite{Mallat2000}. These last relations are the starting point to develp next our main results.
\section{Main results}
In the present paper, our aim is to introduce 'new' wavelet functions (father and mother wavelets) and associatd multiresolution analysis using the well known Farey map. We will show that such a map permits (as in the case of Haar, Morlet, Gaussian functions) to develop a multiresolution analysis on $\mathbb{R}$ and to conduct important applications and algorithms in applied contexts.

\begin{figure}[http]
\begin{center}
\includegraphics[scale=0.75]{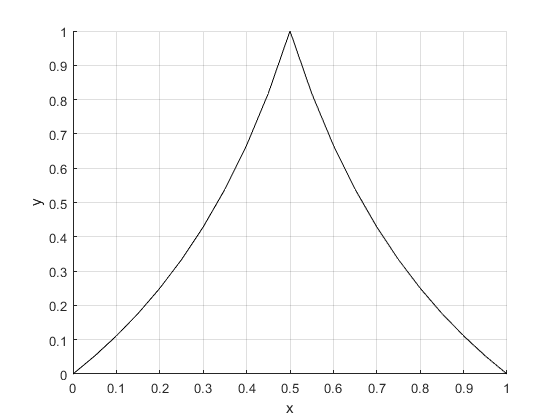}
\caption{The Farey map $F$.}\label{FareyMapF}
\end{center}
\end{figure}

The Faray map (slightly modified) is defined by 
$$
F(x)=\left\{\begin{array}{lll}
\displaystyle\frac{1}{4\log2-2}\displaystyle\frac{x}{1-x}&,&0\leq x\leq\frac{1}{2},\\
\displaystyle\frac{1}{4\log2-2}\displaystyle\frac{1-x}{x}&,&\frac{1}{2}\leq x\leq1.
\end{array}\right.
$$
We may define in a general way a generalized Farey map relatively to a gauge function $h$ instead of $\displaystyle\frac{1}{4\log2-2}\displaystyle\frac{x}{1-x}$ by considering
$$
F_h(x)=\left\{\begin{array}{lll}
h(x)&,&0\leq x\leq\frac{1}{2},\\
h(1-x)&,&\frac{1}{2}\leq x\leq1.
\end{array}\right.
$$
or more generally
\begin{equation}\label{GeneralizedFareyMaphhbar}
F_{h,\overline{h},a}(x)=\left\{\begin{array}{lll}
h(x)&,&0\leq x\leq a,\\
\overline{h}(x)&,&a\leq x\leq1,
\end{array}\right.
\end{equation}
relatively to some suitable functions $h$ and $\overline{h}$ and real number $a$. Figure \ref{FareyMapF} illustrates the Farey map $F$.

In this part we serve of the Farey map above to construct new wavelet function on $\mathbb{R}$ (which may be obviously extended to $\mathbb{R}^n$, $N\geq2$). To do this we consider a bit modification of the Farey map above to obtain the modified Farey map which will be denoted by $\varphi$ and defined as follows
\begin{equation}\label{ModifiedFareyMap}
\varphi(x)=F(\dfrac{1+x}{2})=\left\{\begin{array}{lll}
\displaystyle\frac{1}{(4\log2-2)}\displaystyle\frac{1+x}{1-x}&,&-1\leq x\leq0,\\
\\
\displaystyle\frac{1}{(4\log2-2)}\displaystyle\frac{1-x}{1+x}&,&0\leq x\leq1.
\end{array}\right.
\end{equation}
Related to wavelet theory, the function $\varphi$ is a wavelet-copy-like (or scaling-function-copy-like) of the Farey map $F$ as in wavelet theory the copies of the mother wavelet and the scaling function are always obtained by means of translation/dilation copies. In our case we remark easily that $\varphi=\sqrt2F_{1,-1}$ where for $j,k\in\mathbb{Z}$, $F_{j,k}(x)=\dfrac{1}{2^{j/2}}F(\dfrac{x-k}{2^j})$.

In this section we propose to establish some properties of the modified Farey function $\varphi$ introduced in (\ref{ModifiedFareyMap}) and thus prove consequently that it serves to construct father and mother wavelets.
\begin{Lemma}\label{lemme1}
	The function $\varphi$ satisfies $\widehat{\varphi}(0)=1$.
\end{Lemma}
\textbf{Indeed,} simple computations yield that
$$
\widehat{\varphi}(0)=\dfrac{1}{2\log2-1}\displaystyle\int_{1}^{2}\dfrac{2-t}{t}dt=\dfrac{1}{2\log2-1}\displaystyle\int_{1}^{2}(\dfrac{2}{t}-1)dt=1.
$$
The first point in the construction of wavelet mothers and fathers and consequently associated multi-resolution analysis is to check the existence of Riesz basis $\varphi_k$ (for the eventual space $V_0$), the admissibility, vanishing moments and localization for the mother wavelet $\psi$ and next, the famous 2-scale relation. In this direction, we have a series of results.
\begin{Lemma}\label{lemma-2-scale}
	The function $\varphi$ satisfies the 2-scale relation
	$$
	\varphi(x)=\sqrt2\displaystyle\sum_{k\in\mathbb{Z}}h_k\varphi(2x-k),
	$$
	where $h_1=h_{-1}=\dfrac{1}{3\sqrt2}$, $h_0=\dfrac{1}{\sqrt2}$ and 0 else.
\end{Lemma}
\textbf{Proof.} Because of the disjointness of the supports of the functions $\varphi(.-k)$, $k\in\mathbb{Z}$ appearing in the 2-scale relation we immediately observe that
$$
h_k=0,\,\forall k,\,|k|\geq2.
$$
So the relation in Lemma \ref{lemma-2-scale} above reads
$$
\varphi(x)=\sqrt2h_{-1}\varphi(2x+1)+\sqrt2h_{0}\varphi(2x)+\sqrt2h_{1}\varphi(2x-1).
$$
Next, evaluating the last equation for suitable values of $x$ we obtain
$$
h_1=h_{-1}=\dfrac{1}{3\sqrt2}\quad\mbox{and}\quad h_0=\dfrac{1}{\sqrt2}.
$$
\begin{Corollary}
	For all $\xi\in\mathbb{R}$, we have
	$$
	\widehat{\varphi}(\xi)=\mathcal{M}_0(\dfrac{\xi}{2})\widehat{\varphi}(\dfrac{\xi}{2}),
	$$
	where $\mathcal{M}_0(\xi)=\dfrac{1}{2}\left(1+\dfrac{2}{3}\cos\xi\right)$.
\end{Corollary}
\textbf{Proof.} By applying the Fourier transform to the 2-scale relation in Lemma \ref{lemma-2-scale} above we obtain
$$
\widehat{\varphi}(\xi)=\dfrac{1}{2}\left(\dfrac{1}{3}e^{i\xi/2}+1+\dfrac{1}{3}e^{-i\xi/2}\right)\widehat{\varphi}(\dfrac{\xi}{2})
$$
which reads as
$$
\widehat{\varphi}(\xi)=\dfrac{1}{2}\left(1+\dfrac{2}{3}\cos(\dfrac{\xi}{2})\right)\widehat{\varphi}(\dfrac{\xi}{2}).
$$
Hence, the desired result follows.
\begin{Theorem}\label{TheoremPsi}
	The function $\psi$ defined by
	$$
	\psi(x)=K_0(\sqrt2\displaystyle\sum_{k\in\mathbb{Z}}g_k\varphi(2x-k)),
	$$
	with $g_k=(-1)^{k-1}h_{1-k}$ and $K_0=\dfrac{4\log2-2}{\sqrt{3-4\log2}}$ satisfies
	\begin{itemize}
		\item[\textbf{i.}] $\widehat{\psi}(\xi)=\mathcal{M}_1(\dfrac{\xi}{2})\widehat{\varphi}(\dfrac{\xi}{2})$, where $\mathcal{M}_1(\xi)=\dfrac{K_0}{6}\left(3-2\cos\xi\right)e^{-i\xi}$.
		\item[\textbf{ii.}] $\|\psi\|_2^2=1$.
	\end{itemize}
\end{Theorem}
\textbf{Proof.} Using Lemma \ref{lemme1}, we get $g_k=0$, $\forall,k\leq-1$ and $k\geq3$. So the relation above reads
$$
\psi(x)=K_0\sqrt2g_{0}\varphi(2x)+\sqrt2g_{1}\varphi(2x-1)+\sqrt2g_{2}\varphi(2x-2)
$$
with
$$
g_0=-\dfrac{1}{3\sqrt2},\;g_1=\dfrac{1}{\sqrt2}\quad\mbox{and}\quad g_2=-\dfrac{1}{3\sqrt2}.
$$
Otherwise,
\begin{equation}\label{psialaiddedevarphi1}
\psi(x)=K_0(-\dfrac{1}{3}\varphi(2x)+\varphi(2x-1)-\dfrac{1}{3}\varphi(2x-2)).
\end{equation}
Denote for simplicity $K=\dfrac{1}{(4\log2-2)\sqrt2}$. Consequently,
$$
\|\psi\|_2^2=K_0^2(\dfrac{22}{9K^2}\displaystyle\int_{1}^{2}\left(\dfrac{2-x}{x}\right)^2dx-\dfrac{4}{3K^2}\displaystyle\int_{0}^{1}\left(\dfrac{x(1-x)}{(1+x)(2-x)}\right)^2dx).
$$
Now, standard computations yield that
$$
\|\psi\|_2^2=1.
$$
\begin{Theorem}
	The function $\psi$ is explicitly expressed by
	$$
	\psi(x)=K_1\left\{\begin{array}{lll}
	\dfrac{1}{3}\dfrac{1+2x}{1-2x}&;&-\dfrac{1}{2}\leq x\leq0,\\
	\\
	\dfrac{1}{3}\dfrac{1-2x}{1+2x}-\dfrac{x}{1-x}&;&0\leq x\leq\dfrac{1}{2},\\
	\\
	\dfrac{x-1}{x}-\dfrac{1}{3}\dfrac{1-2x}{3-2x}&;&\dfrac{1}{2}\leq x\leq1,\\
	\\
	-\dfrac{1}{3}\dfrac{3-2x}{1-2x}&;&1\leq x\leq\dfrac{3}{2}.\\
	\end{array}\right.
	$$
	where $K_1=\dfrac{1}{\sqrt{3-4\log2}}$ is the normalization constant.
\end{Theorem}
\textbf{Proof.} We have from Theorem \ref{TheoremPsi}
$$
\psi(x)=K_0(\displaystyle\frac{1}{3}\varphi(2x)-\varphi(2x-1)+\displaystyle\frac{1}{3}\varphi(2x-2)).
$$
We know proceed by evaluating the right hand side quantity piecewise.
\begin{itemize}
	\item On $[-\displaystyle\frac{1}{2},0]$ we get
\end{itemize}
$\varphi(2x-1)=\varphi(2x-2)=0$ and $\varphi(2x)=\dfrac{1+2x}{1-2x}$. Consequently,
$$
\psi(x)=\dfrac{K_0}{3(4\log2-3)}\dfrac{1+2x}{1-2x}=\dfrac{K_1}{3}\dfrac{1+2x}{1-2x}.
$$
\begin{itemize}
	\item On $[0,\displaystyle\frac{1}{2}]$ we get similarly to the previous case
\end{itemize}
$$
\psi(x)=\dfrac{K_0}{4\log2-3}\Bigl(\dfrac{1}{3}\dfrac{1-2x}{1+2x}-\dfrac{x}{1-x}\bigr)=K_1\Bigl(\dfrac{1}{3}\dfrac{1-2x}{1+2x}-\dfrac{x}{1-x}\bigr).
$$
\begin{itemize}
	\item On $[\displaystyle\frac{1}{2},1]$ we get
\end{itemize}
$$
\psi(x)=\dfrac{K_0}{4\log2-3}\Bigl(-\dfrac{1-x}{x}-\dfrac{1}{3}\dfrac{1-2x}{3-2x}\bigr)=K_1\Bigl(\dfrac{x-1}{x}-\dfrac{1}{3}\dfrac{1-2x}{3-2x}\bigr).
$$
\begin{itemize}
	\item On $[1,\displaystyle\frac{3}{2}]$ we get
\end{itemize}
$$
\psi(x)=\dfrac{K_0}{3(4\log2-2)}\dfrac{3-2x}{2x-1}\bigr)=-\dfrac{K_1}{3}\dfrac{3-2x}{1-2x}.
$$
\begin{Theorem}\label{TheoremPsi1}
	The function $\widetilde{\psi}$ defined by
	$$
	\widetilde{\psi}(x)={\psi}(x)-c\chi_{]-1/2,3/2[}(x)
	$$
	with $c=\dfrac{2\log2-1}{6}$, is admissible with one vanishing moment.
\end{Theorem}
\textbf{Proof.} It remains to show the admissibility and vanishing moments. From equation (\ref{psialaiddedevarphi1}) it suffices to show the admissibility the Farey map $F$. Denote
$$
\mathcal{A}_{F}=\displaystyle\int_{0}^{\infty}\displaystyle\frac{|\widehat{F}(\omega)|^2}{\omega}d\omega,
$$
$$
H(\omega)=\displaystyle\int_{0}^{1/2}h(x)e^{i\omega x}dx
$$
and
$$
G(\omega)=e^{-i\omega}H(\omega).
$$
Straightforward computations yield that
\begin{equation}\label{eq1admissibility1}
\widehat{F}(\omega)=H(-\omega)-G(\omega).
\end{equation}
We will now evaluate $\widehat{F}(\omega)$ near the origin $\omega=0$. By applying the equality above, we get
$$
\widehat{F}'(0)=-iG(0)-2G(0)=i\displaystyle\frac{1-4\log2}{4}.
$$
Consequently, near $0$ we get
$$
\widehat{F}(\omega)=i\displaystyle\frac{1-4\log2}{4}\omega+\omega o(\omega),
$$
where $o(\omega)\rightarrow0$ as $\omega\rightarrow0$. Consequently,
$$
\displaystyle\int_{0}^{\eta}\displaystyle\frac{|\widehat{F}(\omega)|^2}{\omega}d\omega<\infty,\;\forall\,\eta>0.
$$
On the other hand, using equality (\ref{eq1admissibility1}), we get
$$
|\widehat{F}(\omega)|\leq2|G(\omega)|.
$$
Next, integration by parts yields that
$$
G(\omega)=\displaystyle\frac{e^{-i\omega/2}}{i\omega}+\displaystyle\frac{1}{i\omega}\displaystyle\int_{1/2}^{1}\displaystyle\frac{e^{-i\omega t}}{t^2}dt.
$$
Consequently,
$$
|G(\omega)|\leq\displaystyle\frac{2}{|\omega|},
$$
which means that
$$
|\widehat{F}(\omega)|\leq\displaystyle\frac{4}{|\omega|}.
$$
As a result, for $\eta>0$ large enough we get
$$
\displaystyle\int_{\eta}^{\infty}\displaystyle\frac{|\widehat{F}(\omega)|^2}{\omega}d\omega\leq\displaystyle\int_{\eta}^{\infty}\displaystyle\frac{16}{\omega^3}d\omega<\infty.
$$
It holds finally that $\mathcal{A}_{F}<\infty$.

Figures \ref{FarayPhi} and \ref{FarayPsi} below illustrate the functions $\varphi$ and $\psi$ respectively.
\begin{figure}[http]
\begin{center}
\includegraphics[scale=0.75]{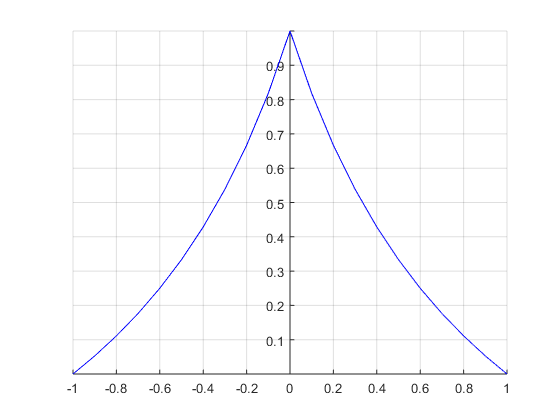}
\caption{Farey scaling function $\varphi$.}\label{FarayPhi}
\end{center}
\end{figure}
\begin{figure}[http]
\begin{center}
\includegraphics[scale=0.75]{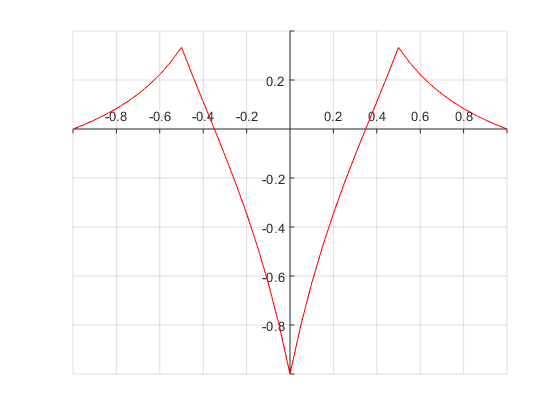}
\caption{Farey wavelet mother $\psi$.}\label{FarayPsi}
\end{center}
\end{figure}
\begin{Remark}
	The supports of the functions $\varphi$ and $\psi$ are compact and satisfy $Support(\varphi)=[N_1,N_2]$ and $Support(\varphi)=[\dfrac{N_1-N_2}{2},\dfrac{N_2-N_1}{2}]$ with $N_1=-1$ and $N_2=1$, which joins the result of Daubechies on compactly supported wavelets.
\end{Remark}
\begin{Definition}
	Let $g\in L^2(\mathbb{R})$. The system of functions $(g_k=g(.-k);\,k\in\mathbb{Z})$ is called a Riesz basis if there exist positive constants $A$ and $B$ such that for any finite set of integers $\Lambda\subset\mathbb{Z}$ and real numbers $\lambda_k$; $k\in\Lambda$, we have
	$$
	A\displaystyle\sum_{k\in\Lambda}\lambda_k^2\leq\|\displaystyle\sum_{k\in\Lambda}\lambda_kg_k\|_2^2\leq B\displaystyle\sum_{k\in\Lambda}\lambda_k^2.$$
\end{Definition}
\begin{Lemma}\label{RieszBasis-phi-k}
	The Farey function $F$ and the modified Farey function $\varphi$ satisfy the following assertions.
	\begin{enumerate}
		\item The system $(F_{k}=F(.-k),\,k\in\mathbb{Z})$ is orthogonal in $L^2(\mathbb{R})$.
		\item The system $(\varphi_{k}=\varphi(.-k),\,k\in\mathbb{Z})$ is only a Riesz system in $L^2(\mathbb{R})$.
	\end{enumerate}
\end{Lemma}
\textbf{Proof.} The first point is a simple consequence of the disjoint supports of the functions $F_{k}$. So, we prove the second. Let $\Lambda\subset\mathbb{Z}$ be a finite set of integers and let $(\lambda_k)_{k\in\Lambda}$ be real numbers. We have
$$
\|\displaystyle\sum_{k\in\Lambda}\lambda_k\varphi_k\|_2^2=
\displaystyle\sum_{l,k\in\Lambda}\lambda_l\lambda_k<\varphi_l,\varphi_k>.
$$
Due to the supports of the functions $\varphi_k$ the last equality reads as
$$
\|\displaystyle\sum_{k\in\Lambda}\lambda_k\varphi_k\|_2^2=\displaystyle\sum_{k\in\Lambda}\left(\lambda_k\lambda_{k-1}<\varphi_k,\varphi_{k-1}>+
\lambda_k^2\|\varphi_k\|_2^2+\lambda_k\lambda_{k+1}<\varphi_k,\varphi_{k+1}>\right).
$$
As $\varphi_k$ are positive functions for all $k$, we obtain
$$
\|\displaystyle\sum_{k\in\Lambda}\lambda_k\varphi_k\|_2^2\geq\displaystyle\sum_{k\in\Lambda}\lambda_k^2\|\varphi_k\|_2^2=\|\varphi\|_2^2\displaystyle\sum_{k\in\Lambda}\lambda_k^2.
$$
On the other hand, using Cauchy-Schwartz inequality we obtain
$$
\displaystyle\sum_{k\in\Lambda}\left(\lambda_k\lambda_{k-1}<\varphi_k,\varphi_{k-1}>+
\lambda_k\lambda_{k+1}<\varphi_k,\varphi_{k+1}>\right)\leq\left(<\varphi,\varphi_{-1}+\varphi_{1}>\right)\displaystyle\sum_{k\in\Lambda}\lambda_k^2.
$$
As a result, by taking $A=\|\varphi\|_2^2$ and $B=\|\varphi\|_2^2+<\varphi,\varphi_{-1}+\varphi_{1}>$ we obtain $0<A<B<\infty$ and
$$
A\displaystyle\sum_{k\in\Lambda}\lambda_k^2\leq\|\displaystyle\sum_{k\in\Lambda}\lambda_k\varphi_k\|_2^2\leq B\displaystyle\sum_{k\in\Lambda}\lambda_k^2.
$$
So as the Lemma.
\begin{Corollary}
	The function $\Gamma_\varphi(\omega)=\displaystyle\sum_{k\in\mathbb{Z}}|\widehat{\varphi}(\omega+2k\pi)|^2$ is bounded on $\mathbb{R}$.
\end{Corollary}
\textbf{Proof.} With the same notations in Lemma \ref{RieszBasis-phi-k} denote
$$
H(x)=\displaystyle\sum_{k\in\Lambda}\lambda_k\varphi_k(x)
$$
and
$$
\widetilde{H}(\xi)=\displaystyle\sum_{k\in\Lambda}\lambda_ke^{-ik\xi}.
$$
We have
$$
\displaystyle\int_{\mathbb{R}}|\displaystyle\sum_{k\in\Lambda}\lambda_ke^{-ik\xi}|^2dx=\displaystyle\int_{\mathbb{R}}|H(x)|^2dx=\dfrac{1}{2\pi}\displaystyle\int_{\mathbb{R}}|\widehat{H}(\xi)|^2d\xi.
$$
Observe now that
$$
\widehat{H}(\xi)=\widetilde{H}(\xi)\widehat{\varphi}(\xi).
$$
Hence,
$$
\displaystyle\int_{\mathbb{R}}|\displaystyle\sum_{k\in\Lambda}\lambda_ke^{-ik\xi}|^2dx=\dfrac{1}{2\pi}\displaystyle\int_{\mathbb{R}}|\widetilde{H}(\xi)|^2|\widehat{\varphi}(\xi)|^2d\xi=\dfrac{1}{2\pi}\displaystyle\sum_{l\in\mathbb{Z}}\displaystyle\int_{2\pi l}^{2\pi(l+1)}|\widetilde{H}(\xi)|^2|\widehat{\varphi}(\xi)|^2d\xi.
$$
Observing next that $\widetilde{H}$ is $2\pi$-periodic we get
$$
\displaystyle\int_{\mathbb{R}}|\displaystyle\sum_{k\in\Lambda}\lambda_ke^{-ik\xi}|^2dx=\dfrac{1}{2\pi}\displaystyle\int_{0}^{2\pi}|\widetilde{H}(\xi)|^2\displaystyle\sum_{l\in\mathbb{Z}}|\widehat{\varphi}(\xi+2\pi l)|^2d\xi=\dfrac{1}{2\pi}\displaystyle\int_{0}^{2\pi}|\widetilde{H}(\xi)|^2\Gamma_\varphi(\xi)d\xi.
$$
Now, as the set $(\varphi_k)_k$ is a Riesz system on $L^2(\mathbb{R})$, we deduce that $\Gamma_\varphi$ is bounded.
\begin{Definition}
	The function $\Gamma_\varphi$ is called the overlap function associated to the system $(\varphi_k)$ or to the function $\varphi$.
\end{Definition}
\begin{Proposition}
	Consider the function $\Phi\in L^2(\mathbb{R})$ defined by its Fourier transform $\widehat{\Phi}=\dfrac{\widehat{\varphi}}{\sqrt{\Gamma_\varphi}}$. Then, the system $(\Phi_k=\Phi(.-k))_{k\in\mathbb{Z}}$ is orthonormal in $L^2(\mathbb{R})$.
\end{Proposition}
\textbf{Proof.} We have for all $l,k\in\mathbb{Z}$,
$$
<\Phi_k,\Phi_l>=\dfrac{1}{2\pi}<\widehat{\Phi}_k,\widehat{\Phi}_l>=
\dfrac{1}{2\pi}\displaystyle\int_{\mathbb{R}}|\widehat{\Phi}(\xi)|^2e^{-i(k-l)\xi}d\xi.
$$
Consequently,
$$
<\Phi_k,\Phi_l>=\dfrac{1}{2\pi}\displaystyle\sum_{n\in\mathbb{Z}}\displaystyle\int_{2n\pi}^{2(n+1)\pi}|\widehat{\Phi}(\xi)|^2e^{-i(k-l)\xi}d\xi=
\dfrac{1}{2\pi}\displaystyle\int_{0}^{2\pi}\displaystyle\sum_{n\in\mathbb{Z}}|\widehat{\Phi}(\xi+2n\pi)|^2e^{-i(k-l)\xi}d\xi.
$$
Observe newt that
$$
\displaystyle\sum_{n\in\mathbb{Z}}|\widehat{\Phi}(\xi+2n\pi)|^2=\dfrac{\displaystyle\sum_{n\in\mathbb{Z}}|\widehat{\varphi}(\xi+2n\pi)|^2}{\Gamma_\varphi(\xi)}=1.
$$
We obtain
$$
<\Phi_k,\Phi_l>=\dfrac{1}{2\pi}\displaystyle\int_{0}^{2\pi}e^{-i(k-l)\xi}d\xi=\delta_{lk}.
$$
\begin{Lemma}
	The Fourier transform of the modified Farey map $\varphi$ is expressed by
	$$
	\widehat{\varphi}(\xi)=4\cos(\xi)Ci(\xi)-4\sin(\xi)Si(\xi)-2\dfrac{\sin\xi}{\xi};\;\forall\,\xi\not=0
	$$
	where $Ci(\xi)=\displaystyle\int_{\xi}^{2\xi}\dfrac{\cos t}{t}dt$ and $Si(\xi)=\displaystyle\int_{\xi}^{2\xi}\dfrac{\sin t}{t}dt$.
\end{Lemma}
\textbf{Proof.} Assume that $\xi>0$. Elementary calculus yield that
$$
\widehat{\varphi}(\xi)=2\displaystyle\int_{1}^{2}\dfrac{2-t}{t}\cos(\xi(t-1))dt,
$$
which may be written as
$$
\widehat{\varphi}(\xi)=4\displaystyle\int_{1}^{2}\dfrac{\cos(\xi(t-1))}{t}dt
-2\displaystyle\int_{1}^{2}\cos(\xi(t-1))dt.
$$
The last integral is an easy form. Applying standard trigonometric rules for the first one we obtain
$$
\widehat{\varphi}(\xi)=4\cos(\xi)\displaystyle\int_{1}^{2}\dfrac{\cos(\xi t)}{t}dt
-4\sin(\xi)\displaystyle\int_{1}^{2}\dfrac{\sin(\xi t)}{t}dt
-2\dfrac{\sin\xi}{\xi}.
$$
Next, taking $u=\xi t$ for the last integrals we obtain
$$
\widehat{\varphi}(\xi)=4\cos(\xi)Ci(\xi)-4\sin(\xi)Si(\xi)-2\dfrac{\sin\xi}{\xi}.
$$
For $\xi<0$, it suffices to see that $\widehat{\varphi}$ is in fact an even function.
\begin{Lemma}
	The modified Farey wavelet mother $\psi$ satisfies
	$$
	\widehat{\psi}(0)=0\;\;\hbox{and}\;\;\widehat{x\psi}(0)=\log2-\displaystyle\frac{3}{4}.
	$$
\end{Lemma}
\textbf{Proof.} Simple calculus yields that
$$
\widehat{\psi}(0)=\displaystyle\int_{0}^{1/2}h(x)dx-\displaystyle\int_{1/2}^{1}h(1-x)dx=0.
$$
Similarly,
$$
\widehat{x\psi}(0)=\displaystyle\int_{0}^{1/2}xh(x)dx-\displaystyle\int_{1/2}^{1}xh(1-x)dx=\displaystyle\int_{1/2}^{1}\displaystyle\frac{(1-2x)(1-x)}{x}dx=\log2-\displaystyle\frac{3}{4}.
$$
\section{Conclusion}
Wavelet theory has known a great success since its discovery. Mathematically, it provides for function spaces good bases allowing their decomposition into spices associated to different horizons known as the levels of decomposition. A wavelet basis is a family of functions obtained from one function known as the mother wavelet, by translations and dilations. This makes the finding of wavelet mothers of great interest. The analysis of a given function using wavelets passes through the so-called wavelet transform or wavelet coefficient. It is a quantity obtained by a convolution product between the function to be analyzed and the copies of the analyzing wavelet mother.

In the present work one motivation was to construct indeed a wavelet mother
starting from an exploitation of the characteristics of the well known Farey map. Well known characteristics in wavelet theory such aa admissibility and vanishing moments rules, compact support have been established for the new wavelet.

Many extension may be addressed as future direction for the present work.
\begin{itemize}
	\item Exploit more the characteristics of the new Farey wavelet such as its continuous wavelet transform, Fourier-Plancherel type rule as well as Parceval formula.
	\item Associate a discrete wavelet transform for the new Farey wavelet.
	\item Construct suitable multi-resolution analysis.
	\item Develop concrete applications of the new wavelet framework to show the utility of the new constructed wavelet mother.
\end{itemize}

\end{document}